\documentclass[12pt]{article}
\usepackage{amscd,amsfonts,amsmath,amssymb,amstext,amsthm}
\title {On the Cycle Spaces Associated to Orbits of Semi-simple Lie Groups}

\author{B. Ntatin\footnote{Dissertation submitted in partial fulfillment for
  the award of a Ph.D. degree in Mathematics at the Ruhr Universit\"at Bochum,
  Germany.}}
\newcommand{\s} {{\smallskip \noindent}}
\newcommand{\m} {{\medskip \noindent}}
\renewcommand{\b} {{\bigskip \noindent}}
\newtheorem{thm} {Theorem} [section]

\newtheorem{lem} [thm]{Lemma}
\newtheorem{prop}[thm]{Proposition}
\newtheorem{cor} [thm]{Corollary}

\theoremstyle{definition}
\newtheorem{defn} {Definition} [section]

\theoremstyle{remark}

\begin{document}
\maketitle
\begin{abstract}
Let $G$ be a semi-simple Lie group and $Q$ a parabolic subgroup of its
complexification $G^\mathbb C$, then $Z=G^\mathbb C/Q$ is a compact complex
$G^\mathbb C$-homogeneous  manifold. The group $G$ as well as $K^\mathbb C$,
the complexification of the maximal compact subgroup of $G$, acts naturally on
$Z$ with finitely many orbits. For any $G$-orbit $\gamma$, there exists a
$K^\mathbb C$-orbit $\kappa$ such that the intersection $\gamma\cap \kappa$
is non-empty and compact. Considering cycle intersection at the boundary of a
$G$-orbit, a definition of the cycle space associated to any $G$-orbit is
given. Using methods involving Schubert varieties and Schubert slices together
with geometric properties of a certain complementary incidence
hypersurface, the cycle space associated to an arbitrary $G$-orbit $\gamma$ is
completely characterised. In particular, it is shown that all the cycle spaces
except in a few Hermitian cases are equivalent to the domain $\Omega_{AG}$. In
the exceptional Hermitian cases, the cycle spaces are equivalent to the
associated bounded domain.
\end{abstract}

\section {Introduction}
If $G$ is a semi-simple Lie group and $Q$ a parabolic subgroup of its
complexification $G^\mathbb C$, then the $G^\mathbb C$-homogeneous, complex
manifold $Z=G^\mathbb C/Q$ is a  projective-algebraic manifold. Moreover, $G$
has only finitely many orbits, thus at least an open orbit in $Z$ (\cite
{W1}. Theorem 2.6).

\m
The complexification $K^\mathbb C$ of the maximal compact subgroup $K\subset
G$  also has finitely many orbits in $Z$. Let $Orb_Z(G)$
(resp. $Orb_Z(K^\mathbb C)$) denote the set of $G$-orbits (resp. $K^\mathbb
C$-orbits) in $Z$. If $\kappa \in Orb_Z(K^\mathbb C)$ and $\gamma \in
Orb_Z(G)$, then $(\gamma, \kappa)$ is said to be a dual pair if the
intersection $\kappa \,\cap \,\gamma$ is non-empty and compact. Duality between
$G$ and $K^\mathbb C$-orbits is proved in (\cite {M}, see
also \cite {BL} and \cite {MUV}):\\
{\it For every $\gamma \in Orb_Z(G)$ there exists a unique
$\kappa \in Orb_Z(K^\mathbb C)$ such that $(\gamma ,\kappa )$
is a dual pair and vice versa.}

\m
Every open $G$-orbit contains a unique $K$-orbit which is a complex manifold
(see, \cite {W1}, Lemma 5.1). This is duality for open orbits. Let $D$ be an
open $G$-orbit in $Z$, define the cycle space associated to $D$ to be the
connected component $\Omega_W(D)$ of the set
$$
\{g(C_0):g\in G^\mathbb C \text { and } g(C_0)\subset D\},
$$
where $C_0$ denotes the base cycle. This set is a
$G$-invariant domain in the $G^\mathbb C$-homogeneous, complex manifold
$\Omega =G^\mathbb C.C_0$ contained in the cycle space ${\cal C}_q(Z)$ of
  $q$-dimensional cycles in $Z$, where $q:=dim_\mathbb c C_0$.

\m
By the procedure of identifying a cycle $g(C_0)\in \Omega_W(D)$ with $g\in
G^\mathbb C$, one can consider that $\Omega_W(D)$ is parametrized by the
group $G^\mathbb C$. In this regard, $\Omega_W(D)$ is then an open subset of
$G^\mathbb C$ which is invariant by the right $K^\mathbb C$-action on
$G^\mathbb C$.

\m
Several authors have been concerned with the problem of
describing  cycle domains  $\Omega_W(D)$ of open $G$-orbits in $Z$, (see for
example, \cite{W1}, \cite{W2}, \cite{WZ}, \cite{HW1}, \cite {HW2} and  \cite
{H} among others). The cycle space $\Omega_W(D)$ associated to an open
$G$-orbit $D$ is well understood. With the exception of a few hermitian cases,
it has been shown (\cite {HW1},\cite {FH}) that the cycle space $\Omega_W(D)$
for any open $G$-orbit $D$ in any flag manifold $Z$ agrees with a
certain domain $\Omega_{AG}$ introduced in (\cite {AG}) independent of $D$ and
$Z$. This domain $\Omega_{AG}$ is an open neighborhood of the Riemannian
symmetric space $G/K$ in $G^\mathbb C/K^\mathbb C$ on which the $G$-action is
proper (\cite {AG}).
 
\m
It is appropriate to consider an analogous definition for the cycle
space of lower-dimensional $G$-orbits in $Z$.  For a dual pair 
$(\gamma,\kappa)\in Orb_Z(G)\times Orb_Z(K^\mathbb C)$, define the {\it {cycle
    space}} $C(\gamma)$ associated to
 $\gamma$ as the connected component containing the identity of the
interior of the set 
\begin {equation*}
\{ g\in G^\mathbb C: g(\kappa )\cap \gamma \text { is non-empty and compact}\}.
\end {equation*}
It is clear that for $\gamma =D$ an open $G$-orbit, $C(\gamma)$ agrees with
$\Omega_W(D)$.

\m
Apriori, however, it is not clear that $C(\gamma)$ is a non-empty  open
subset of $G^\mathbb C$ containing the identity. It could
happen that $g(\text {c}\ell(\kappa))$ intersects the boundary $\text
{bd}(\gamma)$ of $\gamma$ in such a way that an arbitrarily  small perturbation
of $g(\kappa)$ has non-compact intersection with $\gamma$. In Subsection 1.2,
we prove that $C(\gamma)$ so defined above is indeed a non-empty open subset of
$G^\mathbb C$ containing the identity. 

\m
Although $C(\gamma)$ is by definition an open subset of $G^\mathbb C$, one can
think of elements of
$C(\gamma)$ as cycles in ${\cal C}_q(Z)$ in the following way. Since
$C(\gamma)$ is $K^\mathbb C$-invariant, one often regards it generically as
being in the
affine homogenous space $G^\mathbb C/K^\mathbb C$. Furthermore, the $G^\mathbb
C$-action on ${\cal C}_q(Z)$ is algebraic and the $G^\mathbb C$-isotropy
subgroup $G^\mathbb C_{C_0}$ at the base  point $C_0:=\text {c}\ell(\kappa)$
contains $K^\mathbb C$. One can show that the orbit $G^\mathbb C.C_0$ is
either  
$G^\mathbb C/{\bar K}^\mathbb C$, where ${\bar K}^\mathbb C$ is a finite
extension of $K^\mathbb C$, or one of the compact Hermitian symmetric spaces
$X_{\pm}:=G^\mathbb C/P_{\pm}$. Here, $P_{\pm}$ are the unique parabolic
subgroups of $G^\mathbb C$ containing $K^\mathbb C$ and $G^\mathbb C/P_-$
(resp. $G^\mathbb C/P_+$) is the compact Hermitian symmetric space dual to and
containing the Hermitian bounded domain ${\cal B}$ (resp. ${\cal {\bar B}}$) as
  a $G$-orbit.

\m
Our main aim in this work is to give a detailed description of the
cycle space $C(\gamma)$ associated to any $G$-orbit $\gamma$.
We prove the following 
~\\
\begin{thm}\label{main theorem}
If $\text {c}\ell(\kappa)$ is either $P_+$- (or $P_-$)-invariant, then
$C(\gamma)$ is the bounded symmetric domain ${\cal B}$ (or its complex
conjugate ${\cal {\bar B}}$), otherwise $C(\gamma)$ agrees with the domain
$\Omega_{AG}$. 
\end{thm}
~\\
For the proof of this result certain subvarieties are introduced and their
intersections with cycles are studied. A Borel subgroup $B$ of $G^\mathbb C$
containing the factor $AN$ of an Iwasawa decomposition $G=KAN$ of $G$ is
referred to as an Iwasawa-Borel subgroup. The closure $S$ of such a $B$-orbit
${\cal O}$ is called an Iwasawa-Schubert variety. A meromorphic function $f\in
\Gamma(S,{\cal O}(*Y))$, with polar set contained in the variety
$S\setminus {\cal O}$ is then constructed. It is shown that the polar set of
its trace transform ${\cal P}:={\cal P}(Tr(f))$ is a  complex $B$-invariant
hypersurface in the complement of $C(\gamma)$. The polar set ${\cal P}$
consists of a union of hypersurface components of the maximal $B$-invariant
hypersurface $H$ in the complement of $C(\gamma)$. This complementary
$B$-invariant hypersurface $H$ is decisive in the complete characterisation of
$C(\gamma)$.  

\m
It could be possible that the maxmal $B$-invariant hypersurface $H$ in the
complement of $C(\gamma)$ is a lift, i.e., of the form
$H=\pi_+^{-1}(H_+)$ (or $H=\pi_-^{-1}(H_-)$) with respect to the standard
projection $\pi_+:G^\mathbb C/K^\mathbb C\to X_+$ (or $\pi_-:G^\mathbb
C/K^\mathbb C\to X_-$). Here, $H_+$ 
(resp. $H_-$) is the unique $B$-invariant  hypersurface in $X_+$ (resp. $X_-$).

\m
Distinguishing the case where the $B$-invariant hypersurface $H$ is not lift
and the case where $H$ is a lift from either $X_+$ or
from $X_-$, we prove our main results in Section 7 for non-closed
$G$-orbits. The case of the closed $G$-orbit is a special case and is
considered seperately.

\m
If the $B$-invariant hypersurface $H$ is not a lift, then the domain
$\Omega_H$, defined as the connected
component containing the base point in $\Omega:=G^\mathbb C/K^\mathbb C$ of
the set $\Omega\setminus \underset {k\in K}\cup(kH)$, agrees with the
universal domain $\Omega_{AG}$ (\cite {FH}). This
together with certain known results about the cycle spaces of open orbits
(\cite {GM}, \cite {HW1}, \cite {WZ}), leads to the proof of our main
result in this situation.

\m
If the $B$-invariant hypersurface $H$ is a lift, then the base cycle $C_0=\text
{c}\ell(\kappa)$ is $P_+$- (or $P_-$-) invariant. Consequently, the orbit
$G^\mathbb C.C_0$ in the cycle space $C_q(Z)$ is either the symmetric space
$X_+$ or $X_-$. The cycle space $C(\gamma)$
in this case is just a lift of the domain $\Omega_{H_+}$ (or $\Omega_{H_-}$).

\m
The closed $G$-orbit $\gamma_{\text {c}\ell}$ is special in the sense that the
Schubert slices in this case are just points. Consequently, Schubert slice
intersection methods are not helpful. The case of the cycle space
$C(\gamma_{\text {c}\ell})$ of the unique closed $G$ orbit  $\gamma_{\text
  {c}\ell}$ is handled in Section 8. Here, the boundary
of the dual open $K^\mathbb C$-orbit $\kappa_{op}$, is decomposed into
irreducible components $\text {bd}(\kappa_{op})=Z\setminus \kappa_{op}=A_1\cup
\ldots \cup A_k$ and it is shown that in each $A_j$ there is a unique
$K^\mathbb C$-orbit $\kappa_j$ with dual $G$-orbit $\gamma_j$ such that $\text
{c}\ell(\gamma_j)=\gamma_j\dot \cup \gamma_{\text {c}\ell}$. This leads  to
the inclusion $C(\gamma_{\text {c}\ell})\subset C(\gamma_j)$. By showing that
the dual open $K^\mathbb C$-orbit 
$\kappa_{op}$ is neither $P_+$- nor $P_-$-invariant, we obtain the equality
$C(\gamma_{\text {c}\ell})=\Omega_{AG}$ proving the main result for closed
orbits as well.

\m
In conclusion therefore, it is proven that except in a few explicit cases, the
cycle space $C(\gamma)$ associated to any $G$-orbit $\gamma$ in any
$G^\mathbb C$-flag manifold $Z=G^\mathbb C/Q$ agrees with the universal domain
$\Omega_{AG}$. The only exception occurs when the real form $G$ is of
Hermitian type and the base cycle $\text {c}\ell(\kappa)$ is $P_+$- (or
$P_-$-) invariant. Then the cycle space $C(\gamma)$ associated to a nonclosed
$G$-orbit $\gamma$ is the bounded symmetric domain ${\cal B}$ (or its complex
conjugate ${\cal {\bar B}}$).

\m
I heartily thank my surpervisor Prof.\,Dr.\,Dr.h.c. A.T. Huckleberry for
introducing me to this area of research and for always finding time to direct
me in between tight schedules. My thanks also go to Prof.\,Dr. P. Heinzner and
all the members of the department of complex analysis at the Ruhr Universit\"at
Bochum for their support throughout my Ph.D. programe. Finally, I would like
to express my appriciation to the Deutsche Akademische Austauschdienst (DAAD)
for supporting this work.

\section {Preliminaries and Notations}
\subsection{Notations}
Let $G$ be a non-compact semi-simple Lie group which is embedded in its
complexification $G^\mathbb C$ and let $Z=G^\mathbb C/Q$ be a $G^\mathbb
C$-flag manifold, i.e., a compact, homogeneous, algebraic, rational $G^\mathbb
C$-manifold. Here $Q$ is a parabolic subgroup in the sense that it contains a
Borel subgroup. Observe that $G$ acts naturally on every flag manifold
$Z=G^\mathbb C/Q$.  The semi-simple Lie group $G$ decomposes as a finite
direct product of simple groups. This leads to a decomposition of each flag
manifold $Z=Z_1\times...\times Z_k$ as a finite direct product with
irreducible factors $Z_i=G_i^{\mathbb C}/Q_i$, where for each $i$, $Q_i:=Q\cap
G_i$ is a parabolic subgroup of the complexification $G_i^\mathbb C$ of the
simple factors $G_i$. 
Thus a $G$-orbit in $Z$ is a product of $G_i$-orbits
in the corresponding factors $Z_i$. Consequently, in the sequel we will assume
without loss of generality that $G$ is simple. The necessary adjustments for
the semi-simple case are straight-forward.

\m
Fix a Cartan involution $\theta$ of $G$, and extend it as usual
(holomorphically) to $G^\mathbb C$. The fixed point set $K:=G^\theta$ is a
maximal compactly  embedded subgroup  of $G$ and $K^\mathbb C:=(G^\mathbb
C)^\theta$ is its complexification. The compact group $K$ as well as its
complexification $K^\mathbb C$ also act naturally on $Z$ and further,
$G/K$ is a negatively curved Riemannian symmetric space embedded in $G^\mathbb
C/K^\mathbb C$.

\m
Let $Orb_Z(G)$ (resp. $Orb_Z(K^\mathbb C)$) denote the set of $G$-orbits
(resp. $K^\mathbb C$-orbits) in $Z$. It is known that these sets are finite
(\cite {W1}). Since the $G$-action is algebraic, it follows that there is at
least one open $G$-orbit. In fact, the maximal-dimensional $G$-orbits are open
while the minimal-dimensional $G$-orbits are closed. Moreover, there is only
one closed $G$-orbit in $Z$. 
\s

There exists a duality relation between $G$-orbits and $K^\mathbb C$-orbits in
any flag manifold $Z$. We will
use the following  slightly modified version of this duality relationship: {\it
If $\kappa \in Orb_Z(K^\mathbb C)$ and $\gamma \in Orb _Z(G)$, then
$(\kappa ,\gamma )$ is said to be a dual pair if $\kappa \,\cap \,\gamma $
is non-empty and compact}.

\m
If $\gamma $ is an open $G$-orbit, however, then $\kappa $ being dual to
$\gamma $ is equivalent to $\kappa \subset \gamma $. In (\cite {W1})
it is shown that every open $G$-orbit contains a unique compact
$K^\mathbb C$-orbit, i.e., duality at the level of open $G$-orbits.

\m
Duality between $G$- and $K^\mathbb C$-orbits in $Z$ is extended in
(\cite {M}, see also \cite {BL} and \cite {MUV}) to the case of all orbits:
{\it
For every $\gamma \in Orb_Z(G)$ there exists a unique
$\kappa \in Orb_Z(K^\mathbb C)$ such that $(\gamma ,\kappa )$
is a dual pair and vice versa.}

Furthermore, if $(\gamma ,\kappa )$ is
a dual pair then the intersection $\kappa \,\cap \,\gamma $ is transversal
at each of its points and consists of exactly one $K$-orbit.

\m
Let $D$ be an open $G$-orbit in the flag manifold $Z$ and $C_0$ the dual
$K^\mathbb C$-orbit. Since $C_0$ is compact and contained in $D$, it defines a
point in the space of $q$-dimensional compact cycles ${\cal C}_q(D)$, where
$q:=dim_\mathbb CC_0$. By associating $g\in G^\mathbb C$ to the cycle
$g(C_0)$, the connected component $\Omega_W(D)$ of the set  
$$
\{g\in G^\mathbb C :g(C_0)\subset D\}
$$
can be regarded as a family of $q$-dimensional cycles. Of course $\Omega_W(D)$
is invariant by the $K^\mathbb C$ action on $G^\mathbb C$ on the right and
therefore, we often regard it as being in the affine homogenous space
$G^\mathbb C/K^\mathbb C$.

\m
Fix an open $G$-orbit $D$ in $Z=G^\mathbb C/Q$, a base cycle $C_0$ in $D$ and
let $\Omega :=G^\mathbb C.C_0$ denote the corresponding orbit in the cycle
space ${\cal C}_q(D)$. The cycle space $\Omega_W(D)$ associated to open
$G$-orbits in any flag manifold $Z$ has been completely characterized (\cite
{FH}). The following result was proved: 

\begin{thm}(\cite {FH}).
If $\Omega$ is compact, then either $\Omega_W(D)$ consists of a single point
or $G$ is Hermitian and $\Omega_W(D)$ is either the associated bounded
symmetric domin
${\cal B}$ or its complex conjugate ${\cal {\bar B}}$. If $\Omega$ is
non-compact, then regarding $\Omega_W(D)$ as a domain in $G^\mathbb
C/K^\mathbb C$, it follows that $\Omega_W(D)$ agrees with the domain
$\Omega_{AG}$ for every $G$-orbit in every $G^\mathbb C$-flag manifold
$Z=G^\mathbb C/Q$.
\end{thm}

\m
In terms of roots, the domain $\Omega_{AG}$ admits the following
description.
Let ${\mathfrak g}={\mathfrak k}\oplus{\mathfrak p}$ be a
Cartan decomposition of $Lie(G)$, with respect to a compact real form
${\mathfrak g_u}$ of ${\mathfrak g}$. Let ${\mathfrak a}\subset {\mathfrak p}$
be an abelian subalgebra which is maximal with respect to the condition of
being contained in $\mathfrak p$ and $\Phi $ a system of roots on ${\mathfrak
  a}$.  This gives 
rise to an Iwasawa decomposition ${\mathfrak g}={\mathfrak k}+{\mathfrak a}
+{\mathfrak n}$, of ${\mathfrak g}$. For $\alpha$ a root of ${\mathfrak a}$,
let $H_{\alpha}:=\{\xi \in {\mathfrak a}:\alpha(\xi)=\frac {\pi}{2}\}$
and define
$\omega _{AG}$ as the connected component containing $0\in {\mathfrak a}$
of the set which is obtained by removing from ${\mathfrak a}$ the union of all
the affine hyperplanes $H_{\alpha}$ as $\alpha$ runs through the set of
roots. That is, 
 
\begin {equation*}
\omega_{AG}=(\mathfrak a\setminus (\underset{\alpha\in \Phi}\bigcup
H_\alpha))^o=\underset{\alpha \in \Phi}\bigcap \{\xi \in {\mathfrak a}:\vert
\alpha (\xi )\vert < \frac{\pi}{2}\}.  
\end {equation*}

\begin{defn}
Let $\omega_{AG}$ be as above, then the domain $\Omega_{AG}$ is the open
neigborhood of the Riemannian symmetric space $G/K$ in $\Omega:=G^{\mathbb
  c}/K^{\mathbb C}$ given by
$$
\Omega _{AG}=Gexp(i\omega _{AG}).x_0,
$$
where $x_0\in G^\mathbb C/K^\mathbb C$ is the base point. 
\end{defn}

This domain $\Omega_{AG}$ has the property that the $G$-action on it is proper
(\cite {AG}). For a survey of this and other basic properties of
$\Omega _{AG}$ see (\cite {HW2}).

\m
Our main aim in this work is to define and characterize in general the cycle
space associated to any
$G$-orbit $\gamma$ in any $G^\mathbb C$-flag manifold $Z$. With
duality between $G$ and $K^\mathbb C$-orbits in mind, a natural
candidate for the cycle space associated to an arbitrary $G$-orbit $\gamma$
would be the connected component $C(\gamma)$ of the set
\begin {equation*}
\{ g\in G^\mathbb C: g(\kappa )\cap \gamma \text { is non-empty and
  compact}\},
\end {equation*}
for $(\gamma,\kappa)\in Orb_Z(G)\times Orb_Z(K^\mathbb C)$ a dual pair.
This set was introduced in (\cite {GM}) and it was conjectured that the
intersection of all such sets for all $K^\mathbb C$-orbits in all $G^\mathbb
C$-flag manifolds $G^\mathbb C/Q$ agrees with the universal domain
$\Omega_{AG}$. Furthermore, this conjecture was shown (in the same
paper) to be true for classical
and exceptional Hermitian groups by case-by-case considerations.

\m
However, it is not quite clear what happens at the boundary of $\gamma$. Thus
we can not apriori say if the set $C(\gamma)$ is non-empty or not. We will
begin by investigating the intersection of $\text {c}\ell(\kappa)$ wih the
boundary of $\gamma$. This information will lead us to understand boundary
behaviour of {\it cycles} and hence a suitable definition
of the cycle space $C(\gamma)$ that holds for open as
well as for non-open $G$-orbits in any flag manifold $Z$ in general.

\m
Let us now put together some preparatory results. We will closely follow the
notation in (\cite {HW1}). 

\m
Let $B$ denote a Borel subgroup of $G^\mathbb C$ which contains the factor
$AN$ of an Iwasawa decomposition $G=KAN$ of $G$. Such a Borel subgroup is
called an
{\it {``Iwasawa-Borel''}} subgroup of $G^{\mathbb C}$.

\m
Given $B$, an Iwasawa-Borel subgroup of $G^{\mathbb C}$, the closure $S=\text
{c}\ell (O)\subset Z$ of a $B$-orbit $O$ in $Z$ is referred to as an Iwasawa-
Schubert variety. We set $Y:=S\setminus O$. For a fixed Iwasawa-Borel subgroup
$B$, let ${\cal S}$ denote the set of all Schubert varieties and define for
every $\kappa \in Orb_Z(K^\mathbb C)$ the set

\begin {equation*}
{\cal S}_\kappa :=\{ S\in {\cal S}:S\cap \text{c}\ell (\kappa )\not =\emptyset
\}.
\end {equation*}

This set is non-empty since the set of Schubert varieties generate the
integral homology of $Z$.

\m
Given a dual pair $(\gamma $,$\kappa)\in Orb_Z(G)\times Orb_Z(K^\mathbb C)$, a
Borel subgroup $B$ of $G^\mathbb C$ which contains the factor $AN$ of an
Iwasawa-decomposition $G=KAN$ of $G$, a $B$ Schubert variety $S\in {\cal
  S}_\kappa $ and an intersection point $z_0\in \kappa \cap S$, we
refer to $\Sigma=AN.z_0$ as the associated {\it ``Schubert slice''}.

\m
In what follows, fix an Iwasawa decomposition $G=KAN$ and an Iwasawa-Borel
subgroup $B$ containing the factor $AN$. Set $C_0:=\text {c}\ell(\kappa)$, the
closure of the $K^\mathbb C$-orbit $\kappa$, dual to $\gamma \in Orb_Z(G)$.

\begin{lem}
If $p\in \Sigma \cap C_0$, then the tangent space to $\gamma$ at $p$
decomposes as a direct sum
$$
T_p\gamma=T_p(C_0)\oplus T_p(\Sigma).
$$
\end{lem}
\begin{proof}
From the Iwasawa decomposition $G=KAN$, it follows in particular, that
$T_p\gamma=T_p(C_0)+ T_p(AN.p)$. To see that this sum is direct, it just
suffices to count dimensions noting that the orbit $AN.p$ is contained in
$S$.
\end{proof}

\begin{prop}
Let $(\kappa,\gamma )$ be a dual pair and $S\in {\cal S}_\kappa $.
Then
\begin{enumerate}
\item Any Schubert slice  $\Sigma$ is open in $S$
\item At each of its intersection points $S$ is transversal to $\cap C_0$  in
  $\gamma$.
\item  The intersection $S\cap C_0$ is finite and contained in ${\cal
    O}$. Moreover, $S\cap \gamma$ is a finite disjoint union $\dot
    \cup_{j=1}^d \Sigma_j$ of Schubert slices.
\end{enumerate}    
\end{prop}
\begin{proof}
Since $dim_\mathbb C \Sigma + dim_\mathbb C C_0=dim_\mathbb C \gamma$, it
follows
that $dim_\mathbb R S= dim_\mathbb R AN.p$ and as a consequence, $AN.p$ is
open in $S$. Tranversality now follows since we have the direct sum
decomposition $T_p\gamma=T_p(C_0)\oplus T_pS$ of the tangent space.

\m
Any component of $S\cap \gamma$ is $AN$-invariant and since every $AN$-orbit
in $\gamma$ intersects $C_0$, it follows that such an orbit is open in $S\cap
\gamma$. Consequently, such a component is a Schubert slice. It follows that
$S\cap C_0=\{p_1,\ldots,p_d\}$, and $\Sigma_j:=AN.p_j$ are the corresponding
Schubert slices through $p_j$. Hence $S\cap C_0 =S\cap {\cal O}\subset
{\cal O}$ is just the disjoint union of the Schubert slices $\Sigma_j$.

\end{proof}

The following result is implicit in (\cite {HW1}, see Section 5). There it was
only proven that the intersection $\Sigma \cap C_0$ is finite.
\begin{prop}
The intersection $\Sigma \cap C_0$ consists of exactly one point for any
Schubert slice $\Sigma$.
\end{prop}
\begin{proof}
Let $\Sigma$ be a Schubert slice through $p\in \Sigma \cap C_0$. Suppose
$\Sigma$ intersects $C_0$ in another point $p^\prime$. Then 
since $C_0$ is a $K$-orbit, there exists $k\in K$ such that
$k.p^\prime=p$. Since $p\in \Sigma$, there exists $an\in AN$ such that
$(an).p=p^\prime$. It then follows that $kan$ belongs to the $G$-isotropy
subgroup at $p$. 

\m
The map $\alpha:K_p\times (AN)_p\to G_p$ defined by multiplication
$(k,an)\mapsto kan$ is a diffeomorphism (\cite {HW1}) onto a number of
components of $G_p$. However, since $\gamma \,\cap \,\kappa$ is a strong
deformation retract of $\gamma$, it follows that $G_p/K_p$ is connected and
consequently, $\alpha$ is surjective. It therefore follows that $k$
belongs to $K_p$, the $K$-isotropy subgroup at $p$. Thus $k.p^\prime =p^\prime
=p$.
\end{proof}

\subsection {Definition of the cycle Space}
Our aim here is to give a suitable definition of the cycle space associated to
a $G$-orbit in any flag manifold $Z=G^\mathbb C/Q$. We will need the following
result in the sequel.

\m
\begin {prop} \cite {HW1}. \label {basic lemma}
Let $(\kappa,\gamma )$ be a dual pair and $S\in {\cal S}_\kappa $.
Then

\begin {enumerate}

\item $S\cap \text{c}\ell (\kappa )\subset \kappa \,\cap \,\gamma$ .

\item The map $K\times \text {c}\ell(\Sigma)\to \text {c}\ell(\gamma)$, given
  by $(k,z)\mapsto k(z)$, is surjective, that is $K.\text
  {c}\ell(\Sigma)=\text {c}\ell(\gamma)$.

\end {enumerate}
\end {prop}

\m
\begin{cor}\label {contains}
Every $p\in \text {c}\ell(\gamma)$ is contained in some Schubert variety $S\in
{\cal S}_\kappa$.
\end{cor}
\begin{proof} Let $B$ be an Iwasawa-Borel subgroup containing the factor $AN$
of some given Iwasawa decomposition $G=KAN$ of $G$. Let $z_0$ be the base
point of $\gamma$ so that $z_0\in \kappa \cap \gamma$ and $S\in {\cal
S}_\kappa$ be an Iwasawa-Schubert variety through $z_0$, that is , $S=\text
{c}\ell(B.z_0)$. Furthermore, let $\Sigma =AN.z_0$ be a Schubert slice
through $z_0$. Suppose that $p\in \text {c}\ell(\gamma)$, then by
Prop.$\,$\ref{basic lemma}, there exists $k\in K$ such that $p\in k.\text
{c}\ell(\Sigma)$. The $K$-conjugate $kBk^{-1}$ of $B$ contains the conjugate
$kANk^{-1}$ of $AN$ and as a consequence, $p\in k.\text
{c}\ell(kANk^{-1}.z_0)=\text {c}\ell((kAk^{-1})(kNk^{-1})).z_0=\text
{c}\ell(\tilde {A}\tilde {N}.z_0)\subset \text {c}\ell(\tilde {B}.z_0)=\tilde
{S}$. Here $\tilde {S}$ is another Schubert variety which is the closure of
the orbit $\tilde {\cal O}$ of the Iwasawa-Borel subgroup $\tilde {B}$
containing the factor $\tilde {A}\tilde {N}$ of some other Iwasawa
decomposition of $G$.
\end{proof}

\begin{lem}\label {equality}
$\text{c}\ell(\kappa)\cap \text{c}\ell (\gamma )=\kappa \,\cap \,\gamma $.
\end{lem}

\begin{proof}
Of course we have that $\kappa \cap \gamma\subset \text{c}\ell(\kappa)\cap
\text{c}\ell (\gamma )$ and so it is sufficient to prove the opposite
inclusion. Suppose $p\in \text{c}\ell(\kappa)\cap \text{c}\ell (\gamma )$,
then by Cor.$\,$\ref {contains}, there is some $S$ containing the point p, that
is  $p\in \text{c}\ell(\kappa)\cap S$, and this intersection is contained in
$\kappa \cap \gamma$ by the first part of  Prop.$\,$\ref {basic lemma}.

\end{proof}

\begin{lem}\label {neigborhoodlemma}
If $Y$ is any compact set in $Z$ with
$Y\cap \text{c}\ell (\gamma )\subset \gamma $, then there is an
open neighborhood $U$ of the $Id$ in $G^{\mathbb C}$ such that
$g(Y)$ has the same property.
\end{lem}
\begin{proof}
Let $d:Z\times Z\to {\mathbb R}^{\geq 0}$ be any distance function on
$Z$. Since $Y\cap \text{c}\ell (\gamma )\subset \gamma$, and $Y$ is compact,
it follows that the distance between $Y$ and  $\text {bd}(\gamma)$ is
positive, that
is, $d(Y,\text {bd}(\gamma))>0$. Furthermore, the function $\beta:G^{\mathbb
  C}\to {\mathbb R}^{\geq 0}$, given by $g\mapsto d(g(Y),\text {bd}(\gamma))$
is continuous. Since $\beta(Id)>0$, it follows that there exists a neigborhood
$U$ of the identity $Id\in G^{\mathbb C}$ such that $d(gY,\text
{bd}(\gamma))>0$ for $g\in U$. Consequently, $g(Y)\cap \text
{c}\ell(\gamma)\subset \gamma$, for any $g\in U$.

\end{proof}

\begin{cor}\label {interio point}
The identity is an interior point of the set
\begin {equation*}
C:=\{ g\in G^\mathbb C: g(\kappa )\cap \gamma \text { is non-empty and
  compact}\}.
\end {equation*}
\end{cor}
\begin{proof}
We have by Lemma \ref {equality}, that $\text{c}\ell(\kappa )\cap \text{c}\ell
(\gamma )=\kappa \cap \gamma $ which is contained in $\gamma$. Since $\text
{c}\ell(\kappa)$ is compact, Lemma \ref {neigborhoodlemma} implies that there
exists a neigborhood $U$ of the identity in  $G^{\mathbb C}$ such that
$g(\text{c}\ell(\kappa)) \cap \text{c}\ell (\gamma )=g(\kappa) \cap \gamma
\subset \gamma$ for $g\in U$. Since $\gamma$ and $\kappa$ are
dual, the intersection $\gamma\cap \kappa$ is tranversal. Consequently, for
$g$ in a possibly smaller neigborhood as $U$, the intersection $g(\kappa)\cap
\gamma$ remains transversal and
therefore non-empty. It now follows that the intersection
$g(\text{c}\ell(\kappa)) \cap \text{c}\ell (\gamma )$ for $g\in U$ is
non-empty. This implies that the identity is an interior point of the set $C$.
\end{proof}

\begin{defn}
Let $(\gamma,\kappa)\in Orb_Z(G)\times Orb_Z(K^{\mathbb C})$ be a dual
pair. The {\it {cycle space}} $C(\gamma)$ associated to a $G$-orbit $\gamma$ is
the connected component containing the identity of the interior of the set
\begin {equation*}
C=\{ g\in G^\mathbb C: g(\kappa )\cap \gamma \text { is non-empty and
  compact}\}.
\end {equation*}
\end{defn}

It is now clear that $C(\gamma )$ is a  non-empty open subset of $G^\mathbb C$
since the identity belongs to the set $C$. Furthermore, if $\gamma=D$ is an
open $G$-orbit in $Z$, then $C(\gamma)$ agrees with the cycle  domain
$\Omega_W(D)$ introduced in (\cite {WeW}).

\m
Of course elements of $C(\gamma )$ are transformations and not
cycles in the sense of points in the cycle space ${\cal C}_q(Z)$.
There will be occasions, however, where we really want to think of
an element of $C(\gamma )$ as a cycle in the latter sense.

\s
For this we recall that the $G^\mathbb C$-action on ${\cal C}_q(Z)$
is algebraic and therefore the orbit $G^\mathbb C.C_0$  can be identified with
the algebraic homogeneous space $G^\mathbb C/G^\mathbb C_{C_0}$.

\m
It follows that the transformation group variant $C(\gamma )$
is invariant under right-multiplication by $G^\mathbb C_{C_0}$.
Thus, if we wish to think of cycles as being in ${\cal C}_q(Z)$,
we replace $C(\gamma )$ by $C(\gamma )/G^\mathbb C_{C_0}$.

\m
Now, the isotropy subgroup $G^\mathbb C_{C_0}$ always contains $K^\mathbb C$.
In fact, if $G$ is not of Hermitian type, then $K^\mathbb C$
is maximal in $G^\mathbb C$ in the sense that the only proper 
subgroups which contain it are finite extensions $\tilde K^\mathbb C$.
Thus in the nonhermitian case $G^\mathbb C_{C_0}$ is at most
a finite extension of $K^\mathbb C$.

\m
In the Hermition case $K^\mathbb C$ is properly contained in
the parabolic subgroup $P_+$ or $P_-$, where $G^\mathbb C/P_+$ and $G^\mathbb
C/P_-$ are the associated compact Hermitian symmetric spaces. Thus it is quite
possible
that $G^\mathbb C_{C_0}$ is one of these subgroups.  For example,
if $x_+$ is the base point in $G^\mathbb C/P_+$ and $\gamma =G.x_+$ is
its (open) orbit, then $\kappa $ is just the base point and
$G^\mathbb C_{C_0}=P_+$.

\m
In the sequel we will sometimes regard $C(\gamma )$ as being
in ${\cal C}_q(Z)$, i.e., in $G^\mathbb C/G^\mathbb C_{C_0}$, where
$G^\mathbb C_{C_0}$ is one of the groups described above.
On occasion, either by going to a finite cover or pulling
back by one of the fibrations $G^\mathbb C/K^\mathbb C\to G^\mathbb C/P_\pm $
we will regard it in $G^\mathbb C/K^\mathbb C$.

\m
Our work here makes use of results and methods from (\cite {BK}, \cite {FH},
\cite {HW1}, \cite {HW2})  which together with knowledge of the
intersection of cycle domains for the open orbits in $G^\mathbb C/B$
(\cite {GM}), implies the inclusion $\Omega _{AG}\subset C(\gamma )$ for all
$\gamma \in Orb_Z(G)$. This inclusion plays an important role in our proofs.

\section {Characterization of $C(\gamma)$ by Schubert slices}
Throughout this section, $C_0$ will denote $\text {c}\ell(\kappa)$ where
$(\gamma,\kappa)\in Orb_Z(G)\times Orb_Z(K^{\mathbb C})$ is a dual
pair. Futhermore, we assume that the $G$-orbit $\gamma$ under
consideration is not closed. We give a characterization of the cycle space
$C(\gamma)$ associated to $\gamma$ by means of cycle intersection
with Schubert slices. Our goal here is to prove the following

\begin{prop}\label {main prop}
Let $\Sigma $ be a Schubert slice and suppose that $\{g_n\}$ is a sequence in
$C(\gamma)$ such that
$g_n\rightarrow g$ with $g_n(\text {c}\ell(\kappa))\cap \Sigma=
\{p_n\}$ and $p_n$
diverges in $\Sigma$. Then $g\not \in C(\gamma)$.
\end{prop}

In order to prove the above result, we need some preparations.

\begin {lem} \label {int.number}
For all $g\in C(\gamma )$ the number of
points in the intersection $g(\kappa )\cap \Sigma $
is bounded by the intersection number $[S].[\text{c}\ell(\kappa )]$.
\end {lem}
\begin {proof} Since $\Sigma $ can be regarded as a domain in
${\cal O}\cong \mathbb C^m$ and $g(\kappa )\cap \gamma $ is compact,
it follows from the maximum principle that $g(\kappa )\cap \Sigma $
is finite and of course it is then bounded by the intersection
number $[S].[\text{c}\ell (\kappa )]$.
\end {proof}

\m
Define a subset ${\cal I}$ of $G^\mathbb C$ as the connected component
containing the identity of the interior of the set
\begin {equation*}
\{ g\in G^{\mathbb C}:\vert g(\kappa )\cap \Sigma \vert =1
\text{ for all } \Sigma \}.
\end {equation*}
Observe that since $\text {c}\ell(\kappa)\cap \text {c}\ell(\gamma)=\kappa\cap
\gamma$, and $\vert \kappa \cap \Sigma \vert =1$
for all $\Sigma$, it follows that $g\in {\cal I}$ for $g$ sufficiently
close to the identity. Thus ${\cal I}$ is an open subset of $G^\mathbb C$
containing the identity.
\s

In the definition of ${\cal I}$ above,  {\it for all \ $\Sigma $} means for
all choices of the maximal compact group $K$ and all Iwasawa factors $AN$,
i.e., all Schubert slices which arise by $G$-conjugation of those $\Sigma $
which are connected components of $S \cap \gamma $ for a fixed $S\in {\cal
  S}_\kappa $ 

\m
\begin{lem} \label{transversality}
For all $g\in {\cal I}$, the intersection $g(\kappa)\cap \gamma$ is
transversal.
\end{lem}
\begin {proof}
Let $d$ denote the intersection
number $[S].[\text {c}\ell(\kappa)]$, then the base cycle $\text
{c}\ell(\kappa)$ intersects ${\cal O}$ in exactly $d$ points. Furthermore,
${\cal O}\cap \gamma$ is a disjoint union of Schubert slices $\Sigma_1,\ldots,
\Sigma_d$ with one intersection point in each slice. Thus if $g\in {\cal I}$,
then $g(\text {c}\ell(\kappa))$ intersects each $\Sigma_i$ in exactly one
point as well. If any of such intersection  point were not transversal, then
the intersection number would be too big.
\end{proof}

The following is a consequence of the above Lemma.
\begin {cor} \label {connected}
The intersection $M_g=g(\kappa )\cap \gamma $ is a connected compact manifold
for all $g\in {\cal I}$.
\end {cor}

\m

\begin {prop} \label {empty}
$\text {bd}({\cal I})\cap C(\gamma)=\emptyset$.
\end {prop}
\begin{proof}
Assume by contradiction that $g\in \text {bd}({\cal I)}\cap C(\gamma)$. Let
$\{g_n\}$ be a sequence in $I\cap C(\gamma)$ with $g_n\to g$. Then by
Cor.$\,$\ref{connected}, $M_n:=g_n(\kappa)\cap \gamma$ is a sequence of compact
connected manifolds in $\gamma$.  

\m
Let $\tilde {M}$ denote the limiting set, $\tilde {M}:=\lim M_n$. It follows
that  $\tilde {M}$ is a connected closed subset of $\text {c}\ell(\gamma)$.

\s
Now,
$$\tilde {M}\subset g(\text {c}\ell(\kappa)) \cap \text
{c}\ell(\gamma)=:A\,\dot \cup \,E,
$$
where
$$
A:=(g(\text {bd}(\kappa))\cap \text {c}\ell(\gamma))\cup
(g(\kappa)\cap\text {bd}(\gamma)):=A_1\,\dot \cup \,A_2
$$
and
$$
E:=g(\kappa)\cap \gamma.
$$

The set $A$ is closed, because $A_1$ is the intersection of two closed sets,
and a sequence in $A_2$ which converges in $Z$ will either converge to a
point of $A_2$ or $A_1$.

\m
Since we have assumed that $g\in C(\gamma)$, it follows that $E$ is
compact. Thus
$$
\tilde {M}=A\,\dot \cup \,E
$$
is a decomposition of $\tilde {M}$ into disjoint open subsets of $\tilde
{M}$. Since $\tilde {M}$ is connected and $\tilde {M}\cap A\neq \emptyset$, we
conclude that $\tilde {M}\subset A$.

\m

Consequently, for every relatively compact open neighborhood $U$ of a point
$p\in \gamma$, there exists a positive integer $N=N(U)$ such that
$g_n(\kappa)\cap U =\emptyset$ for all $n>N$.

\m

Now we have assumed that $g\in C(\gamma)$, in particular that $E$ is
nonempty. Hence, for $p\in E$ we can consider an Iwasawa-Schubert variety
$S={\cal O}\,\dot \cup \,Y$ with $p\in {\cal O}$. Since $E$ is compact, the
complex analytic set $g(\text {c}\ell(\kappa))\cap S$ must contain $p$ as an
isolated point. Thus, for $\Sigma $ a Schubert slice through $p$,
 the intersection $g(\kappa)\cap \Sigma$ is isolated at $p$ and
as a consequence,
$g_n(\kappa)$ must have nonempty intersection with any open neighborhood
$U=U(p)$ of $p$ if $n$ is sufficiently large. This is contrary to the above
statement and therefore $g\not \in C(\gamma)$.

 \end{proof}

\begin{prop}\label {containment}
$C(\gamma)\subset {\cal I}$.
\end{prop}
\begin{proof}
By the openess of ${\cal I}$ and the fact that $\text {bd}({\cal I})\cap
C(\gamma)=\emptyset$ (see, Prop.$\,$\ref {empty} above), it follows in
particular that ${\cal I}\subset C(\gamma)$.
\end{proof}

\m

{\it Proof of Proposition \ref {main prop}}. Again we argue by
contradiction. Suppose that $g\in C(\gamma)$. Then by Prop.$\,$\ref
{containment}, $g(\text {c}\ell(\kappa))\cap \Sigma$ consists of exactly one
point, say $q$. Since $p_n$ diverges in $\Sigma$, we may assume that $p_n\to
p\in \text {c}\ell(\Sigma)\setminus \Sigma$.

\m
Now by definition $C(\gamma)$ is open. Thus there exists a small $h\in
G^{\mathbb C}$ with $hg$ still in $C(\gamma)$ and $hg(\kappa)\cap \Sigma$
containing points near $p$ and $q$. Thus, $\vert hg(\kappa)\cap \Sigma \vert
\geq 2$, in violation of $C(\gamma)\subset {\cal I}$.
\hfill $\square $

\begin{cor}\label{sequences in S}
If $\{g_n\}$ is a sequence of cycles in $C(\gamma)$ such that
$g_n\rightarrow g$ and $g_n(\text {c}\ell(\kappa))\cap (\gamma \cap S)$
contains a sequence $p_n$ which diverges in ${\cal O}$, then $g\not \in
C(\gamma)$.
\end{cor}
\begin{proof}
Since ${\cal O}\cap \gamma$ is a finite union of Schubert slices $\Sigma_1\cup
\ldots \cup \Sigma_k$ say, and the sequence $\{p_n\}$ diverges in ${\cal O}$,
it follows that some Schubert slice contains infinitely many points of the
sequence $\{p_n\}$. We may therefore assume that the sequence $\{p_n\}$ is
contained in some fixed Schubert slice $\Sigma$. This implies that
$p_n\rightarrow p\in \text {c}\ell(\Sigma)\setminus \Sigma$, and it follows
from Prop.\,\ref {main prop} that $g\not \in C(\gamma)$.
\end{proof}

\section {The trace transform}

As usual let $Z$ denote an arbitrary flag manifold and
$z_0$ the base point of a non-closed $G$-orbit $\gamma$.
In this section we regard a cycle as being in ${\cal C}_q(Z)$.

\m
For $(\gamma,\kappa)\in \text {Orb}_Z(G)\times \text {Orb}_Z(K^{\mathbb C})$ a
dual pair, let $S\in {\cal }S_{\kappa}$ be an Iwasawa Schubert variety
throught the base point $z_0$. We recall that $S$ is the closure of ${\cal
  O}$, the orbit of an Iwasawa-Borel subgroup $B$. That is, $S$ is the
disjoint union $S={\cal O}\, \dot \cup \,Y$.
\s

Let $\Psi$ consist of all pairs $(S,Y)$ which occur as above, then we refer to
an element $\psi\in \Psi$ as datum for a {\it trace-transform} ${\cal
T}_{\psi}$ (see definition below), with respect to $\psi$.

Define a subset of the space ${\cal C}_q(Z)$ of $q$-dimensional cycles in $Z$
by

$$
\Omega_{\psi}:=\{C\in {\cal C}_q(Z):C\cap Y=\emptyset \}.
$$ 
Since the condition $C\cap Y\neq \emptyset$ defines a closed analytic
subset in ${\cal C}_q(Z)$, it follows that
$\Omega_\psi$ is Zariski open and dense in the irreducible component of ${\cal
  C}_q(Z)$ which contains the base cycle.  

\m

Let ${\cal M}(S)$ denote the ring of meromorphic functions on $S$ and
 $$
\Gamma(S,{\cal O}(*Y)):=\{f\in {\cal M}(S):{\cal P}(f)\subset Y\}
$$
the subring of functions with polar set  ${\cal P}$ contained in $Y$.
\begin {defn}
The trace-transfomation, ${\cal T}_\psi$, with
respect to $\psi \in \Psi $ is defined by
$$
{\cal T}_{\psi}:\Gamma(S,{\cal O}(*Y))\to {\cal M}(\Omega_\psi),~f\mapsto
{\cal T}_{\psi}(f):=\sum_{p\in C\cap S}f(p). 
$$
\end {defn} 
When there is no ambiguity, we will just write ${\cal T}$ for ${\cal
  T}_{\psi}$. 

\m
Observe that regarded as a function, ${\cal T}_{\psi}(f):\Omega_{\psi}\to
  \mathbb C$, the trace-transformation ${\cal T}_{\psi}(f)$ is infact a
composition of two maps. Firstly, the intersection map $I:\Omega_{\psi} \to
  {\cal C}^0_q(Z)$, which associates to a cycle its intersection with $S$, and
is holomorpic. Secondly, the trace map $\text {trace}(f):Sym_k(S\setminus Y)\to
\mathbb C$ given by
$$
\text {trace}(f)(z_1,\ldots,z_k)=\sum^k_{j=i}f(z_j).
$$
which is also holomorphic.

We need the following
\begin {prop}\cite {HSB}. 
Let $\Gamma$ be a closed irreducible subspace of ${\cal C}_q(Z)$ such that
$\Gamma\cap \Omega_{\psi}\neq \emptyset$. In particular,
$\Gamma\cap\Omega_{\psi}$ is a dense, Zariski open set in $\Gamma$. Then
${\cal T}_{\psi}(f)$ is holomorphic on $\Gamma\cap \Omega_{\psi}$ and extends
meromorphically to $\Gamma$.
\end {prop}
\begin {proof}
We have already seen above that ${\cal T}_{\psi}(f)$ is a composition of two
maps, the intersection map and the trace map. Infact, the intersection map
$$
I:\Gamma\cap \Omega_{\psi}\to \text {Sym}^k(S\setminus Y)
$$
is give by $C\mapsto I(C)=C\cap S$ (of course in the cycle sense), where $k\in
\mathbb N$ depends on $\Gamma$. Since the projection map $\Gamma\times
Sym_k(S)\times S\to \Gamma$ is proper and the set
$$
\{(C,(z_1,\ldots ,z_k),z)\in \Gamma\times Sym_k(S)\times S:z\in C
\text { and } z=z_i \text { for some }i\}
$$
is a closed analytic set in $\Gamma\times Sym_k(S)\times S$, we have
that $I$ extends to a meromorphic map $\tilde {I}:\Gamma\to Sym_k(T)$.
The trace map on the other hand, $trace(f)$ map;
$$
trace(f):Sym^k(S\setminus Y)\to \mathbb C
$$
given by $trace(f)(z_i,\ldots ,z_k)=\sum^k_{j=1}f(z_j)$, is
holormorphic and extends as a meromorphic function $\tilde {f}:\text
{Sym}^k(S)\to \mathbb C$. Indeed, in the neigborhood of a point $(z^0_1,
\ldots, z^0_k)\in Sym_k(S)$, where some of the $z^0_j$ belong to $Y$, choose
$g$ holomorphic near $\{z^0_1\}\cup\ldots \cup \{z^0_k\}$ such that $g\equiv
1$ if $z^0_i\notin Y$, then $g\times f$ is holomorphic. Set
$G(z^0_1,\ldots,z^0_k)=\Pi^k_{j=1}g(z_j)$, if $(z_1,\ldots ,z_k)$ is near
$(z^0_1,\ldots, z^0_k)$, then

$$
(G\times \text {trace}(f)(z_1,\ldots, z_k)=\sum^k_{j=1}g(z_1).\dots \widehat
{g(z_i)}.\ldots .g(z_k).(g\times f)(z_k),
$$
(where $\widehat {g(z_i)}$ means omit this factor in the above sumation), is
meromorphic on $Sym_k(S)$ near $(z^0_1,\dots ,z^0_k)$. 

\s
Consequently, $\text {trace}(f)$ extends meromorphically on $Sym_k(S)$
and so ${\cal T}_{\psi}(f)$ is also meromorphic as a composition
$trace(f)\circ I$ of two meromorphic maps on all of $\Gamma$.
\end {proof}

\begin{cor}
The trace-transform ${\cal T}_{\psi}(f)$ extends meromorphically to the
component of the space of $q$-dimensional cycles ${\cal C}_q(Z)$ which
contains $\Omega$.
\end{cor}

\section{Embedding of the Schubert variety}

Let $S\in {\cal S}_\kappa$ be a Schubert variety. It is our aim in this
section to embed $S={\cal O}\,\dot\cup \,Y$  into a projective space so that
$Y$ is the intersection of $S$ with the hyperplane at infinity.

\m

Let $L\to Z$ be any very ample line bundle on $Z$ and let $L|_S\to S$ be its
restriction to S. Without loss of generality, we may assume that $G^\mathbb C$
is simply connected. Then the line bundle $L\to Z$ is $G^{\mathbb
  C}$-homogeneous and consequently, the
restriction to $S$ is $B$-homogeneous as well. Let $\Gamma(L,Z)$,
(resp. $\Gamma(L|_S,S)$) denote the finite dimensional vector spaces of
holomorphic sections of the respective bundles. Then there exists a natural
$B$-equivariant restriction map

$$
r:\Gamma(L,Z)\to \Gamma(L|_S,S) \text { given by } \sigma \mapsto \sigma|_S.
$$
 Since $L\to Z$ is very ample, it yields a holomorphic $B$-equivariant
 embedding onto its image;

$$
\varphi :S\to \mathbb P(Im(r)^*) \text { given by } s\mapsto H_s:=\{\sigma \in
Im(r):\sigma(s)=0\}.
$$
Moerover, if $((\sigma_0,\ldots, \sigma_m))$ is a basis for $Im(r)$, then
$\varphi$ may be defined in coordinates by the following map
$$
S\mapsto \mathbb P_m(\mathbb C), s\mapsto \varphi(s)=[\sigma_0(s),\ldots,
\sigma_m(s)]. 
$$
By the Borel fixed point theorem, $B$ has an eigenvector $r(\sigma_B)$ in
$Im(r)\setminus \{0\}$. Suppose $\sigma_0:=r(\sigma_B)$ is a $B$-eigenvector
and let $H:=\{\sigma_0=0\}$, then $H$ is $B$-invariant.
Using the above notation, we first note the following
\begin{prop}\label {H equal to Y}
 $H\subset Y$.
\end{prop}
\begin{proof}
Since $H$ is $B$-invariant, it follows that $S\setminus H$ is also
$B$-invariant. Since both ${\cal O}$ and $S\setminus H$ are Zariski dense,
their intersection is nonempty. Thus for $s\in {\cal O}\cap
(S\setminus H)$, it follows that the $B$-orbit ${\cal O}=B.s$ is contained in
$S\setminus H$. This implies that $H$ is contained in the complement of ${\cal
 O}$. 
\end{proof}

With respect to the above
embedding $\varphi $ therefore, the open $B$-orbit ${\cal O}$ is embedded in
$\mathbb C^m$ as follows; 
$$
{\cal O}\to \mathbb C^m ,~ s\mapsto
(\frac{\sigma_1}{\sigma_0}(s), \ldots, \frac{\sigma_m}{\sigma_0}(s)),
$$
and as such, it could be possible that some components of $Y$ are not in $H$.
However, we proceed to show that infact $H=Y$. For this, we need the following

\begin{lem}\label{unipotent subgroup}
If $U$ is a unipotent group acting algebraically as a group of affine
transformations on $\mathbb C^n$, then every orbit is closed.
\end{lem}
\begin{proof}
Let $(w_1,\ldots,w_m)$ represent coordinates for $\mathbb C^m$. Assume
that the action is in upper triangular form and that the coordinates are
arranged so that the projection $\pi:\mathbb C^m\to \mathbb C^{m-1}$ given by
$(w_1,\ldots,w_m)\mapsto (w_2,\ldots,w_{m})$ is $U$-equivariant.

\m
Our aim is to show that for $p\in \mathbb C^m$, the orbit $U.p$ is closed. Let
$q:=\pi(p)$ and assume inductively that $U.q$ is closed in $\mathbb
C^{m-1}$. Since $U.q$ is closed, it follows that 
$\pi(\text {c}\ell(U.p))=U.q$.

\m
Suppose firstly that $U(p)$ and $U(q)$ have the same dimension. Then if $U(p)$
were not closed, it would have an orbit of lower dimension on its boundary
which would be smaller than $U(q)$. But this is impossible, because it would
be mapped onto $U(q)$.

\m
On the other hand, suppose that $U(p)$ is dimension-theoretically larger
than $U(q)$, that is one dimension bigger. Now every fiber of $U(p)\to U(q)$
is a copy of $\mathbb C$ which is open in the $\pi$-fiber. If $U(p)$ were not
closed, then its closure in the $\pi$-fiber would contain an additional
point. Consequently, this closure would be $\mathbb P_1$. But the $\pi$-fiber
is a copy of $\mathbb C$ which certainly does not contain a $\mathbb
P_1$. Hence every fiber of $U(p)\to U(q)$ is a $\pi$-fiber and as a
consequence, $U(p)$ is the $\pi$-preimage of the closed set $U(q)$ and thus is
also closed.

\m

Since every $U$-orbit on $\mathbb C$ is clearly closed, we conclude by the
induction hypothesis that all $U$-orbits in $\mathbb C^m$ are closed.
\end{proof}

Now let $U$ denote the unipotent radical of the Iwasawa-Borel subgroup $B$. We
put all the above results together in the following

\begin{prop}\label{embedding}
Let $S\in {\cal S}_{\kappa}$ be a Schubert variety and $L\to Z$ any very ample
line bundle on $Z$ and suppose 
$$
r:\Gamma(L,Z)\to \Gamma(L|_S,S) \text { given by } \sigma \mapsto \sigma|_S
$$
denotes the canonical restriction map. Furthermore, let the embedding of $S$
in to a projective space be given in coordinates by
$$
S\mapsto \mathbb P_m(\mathbb C), s\mapsto [\sigma_0(s),\ldots, \sigma_m(s)]
$$
where $((\sigma_0,\ldots, \sigma_m))$ denotes a basis for $Im(r)$. If
$\sigma_0$ is chosen to be a $B$-eigenvector, then the $B$-invariant
hyperplane $\{\sigma =0\}$ is the complement $Y$ of the $B$-orbit ${\cal
  O}$. 
\end{prop}
\begin{proof}
Since $U$ is the unipotent radical of the Iwasawa-Borel subgroup $B$, it acts
transitively on the $B$-orbit ${\cal O}$. Now the line bundle $L\to Z$ is
very ample and defines a $B$-equivariant embedding $\varphi:S\to \mathbb
P(Im(r)^*)$ such that the zero set of $\sigma_0$, that is $H=\{\sigma_0=0\}$,
is contained in the complement of the open $B$-orbit ${\cal O}$. It follows
that
$$
S\setminus \{\sigma_0=0\}\to \mathbb C^m \text { given by }s\mapsto
(\frac{\sigma_1}{\sigma_0}(s), \ldots, \frac{\sigma_m}{\sigma_0}(s))
$$
is a $U$-equivariant embedding of $S\setminus H$ into $\mathbb C^m$. By Lemma
\ref{unipotent subgroup}, it follows that the $U$-orbit in $S\setminus H$ is
both open and closed in $\mathbb C^m$. Consequently, the $U$-action on
$S\setminus H$ is transitive and therefore we conclude that $H=Y$.
\end{proof}

\m
Since the $B$-orbit ${\cal O}$ is algebraically isomorphic to $\mathbb
C^{m({\cal O})}$, we will later on apply the above result to show that any
  sequence in ${\cal O}$ that converges to a point in $Y$, converges to
  infinity in $\mathbb C^{m({\cal O})}$. More precisely, we have the following
\begin{cor}\label{hyperplane at infinity}
If a sequence $\{p_n\}\in {\cal O}$ converges to a point $p\in Y$ then with
respect to the above embedding, it converges to infinity in $\mathbb C^m$.
\end{cor}
\begin{proof}
Let $[z_0:\ldots :z_m]$ represent homogeneous coordinates for $\mathbb
P_m(\mathbb C)$ such that the unipotent group $U$ fixes the coordinate
$z_0$. The affine action of $U$ is given by the restriction of its
action on $\mathbb P_m(\mathbb C)$ to the complement of $\{z_0=0\}$ with
affine coordinates $(w_1,\ldots,w_m)$ where $w_j:=\frac {z_j}{z_0}$ for
$j=1,\dots,m$. Since the embedding of $S\setminus 
H$ into $\mathbb C^m$ is $U$-equivariant and $H=Y$ by Prop.\,\ref{embedding},
it follows that if $p_n\in {\cal O}$ converges to a point $p\in Y$ then it
converges to $\infty$ in $\mathbb C^m$.
\end{proof}

\section {Incidence variety $I_Y$}

As usual, let $(\gamma, \kappa)\in Orb_Z(G)\times Orb_Z(K^{\mathbb C})$ be a
dual pair and regard the cycle space $C(\gamma)$ as being contained in the
group $G^{\mathbb C}$. In what follows, $p_0\in \kappa\cap \gamma$ will denote
the base point.

\s
For an Iwasawa-Schubert variety $S\in {\cal S}_\kappa$ defined by an
Iwasawa-Borel subgroup $B$, we have the decomposition
$S={\cal O}\, \dot \cup\,Y$ where ${\cal O}$ is the $B$-orbit, and define the
incidence variety

$$
I_Y:=\{g\in G^{\mathbb C}:g(\text {c}\ell(\kappa))\cap Y\neq \emptyset\}.
$$
Observe that $I_Y$ is a complex analytic subset of $G^{\mathbb C}$. Recall that
for $g\notin I_Y$ the intersection $g(\text {c}\ell(\kappa))\cap S$ is finite.

\m

Let $U_S$ denote the subset of $G^{\mathbb C}$ defined by
$$
U_S=\{g\in G^{\mathbb C}:g(p_0)\in S\}
$$
and let $Q$ be the $G^{\mathbb C}$-isotropy subgroup at $z_0\in Z$. Then the
evaluation map
$$
U_S\to S, \text { given by } g\mapsto g(p_0),
$$
is a $Q$-principal bundle. Indeed, it suffice to show that $Q$ acts
transitively and freely on the fibers. So let $F_s=\{g\in U_S:g(p_0)=s\}$
denote a typical fiber over $s\in S$. Fix $g\in F_s$ then we have the
following identification of $Q$ and $F_s$; $Q\ni h\mapsto g.h\in
F_s$. Consequently, the right $Q$-action on the fibers is transitive and free.

\s
Since $S$ is irreducible, it follows that $U_S$ is an irreducible complex
analytic subset of $G^{\mathbb C}$.

Set ${\cal E}:=U_S\cap I_Y$ and observe that ${\cal E}$ is a proper analytic
subset of $U_S$. Since $U_S$ is irreducible, it follows that ${\cal E}$ is
nowhere dense.
\s

Denote by $D_s$ the subset of $S$ such that the fiber over $s\in S$, is
contained in ${\cal E}$, that is,
$$
D_s:=\{s\in S:F_s\subset {\cal E}\}.
$$
This defines a closed proper complex analytic subset of $S$ which is nowhere
dense as shown in the following

\begin{lem}\label {analytic subset of $S$}
The set $D_S$ is a proper complex analytic subset of $S$.
\end{lem}
\begin{proof}
Observe that the map $\pi:U_S\to S$ is a bundle with connected fibers and
that the set ${\cal E}$ is a closed analytic subset in $U_S$. Let ${\cal
  E}=\cup {\cal E}_i$ be the decomposition of ${\cal E}$ in to irreducible
components. For each $i$, let $A_i$ be the subset of ${\cal E}$ such that the
fiber of $\pi|_{{\cal E}_i}$ at $g\in {\cal E}$ is the same as the $\pi$-fiber
at that point. If $A_i$ is nonempty, then it defines a closed analytic subset
of ${\cal E}$ since
$rank_g(\pi|{{\cal E}_i}):=dim_g{\cal E}_i-dim_g\pi ^{-1}(\pi (g))$
is minimal. Thus $A:=\cup A_i$
is a  closed analytic subset in ${\cal E}$. Now $\pi(A)$ is closed and the
restriction of $\pi$ to each irreducible component of $A$ has constant rank,
therefore it follows that $\pi(A)$ is a finite union of analytic sets and thus
analytic.
\end{proof}

\begin{prop} \label{sequences}
Given a point $g\in U_S$ and a sequence $\{p_n\}\subset S\setminus D_S$
converging to $p=g(p_0)\in Y$, there exists a sequence of transformations
$\{g_n\}\subset U_S\setminus {\cal E}$ with $g_n$ converging to $g$ and
$g_n(p_0)=p_n$.
\end{prop}
\begin{proof}
Let $\{U_n\}$ be a sequence of
open subsets of $U_S$ contracting to $g$, that is, $U_n\subset U_{n+1}$ for all
$n$ and $\bigcap U_n=g$. Since $\pi:U_S\to S$ is an open mapping,
it follows that $V_n:=\pi(U_n)$ 
is a sequence of open neigborhoods of $p$. Consequently, we can renumber the
sequence $p_n$ such that $p_n\in V_n$ for each $n$. Since the set ${\cal E}$ is
a nowhere dense analytic subset of $U_S$ (see, Lemma \ref{analytic subset of
  $S$}) and $p_n\not \in D_S$, it follows that ${\cal E}\cap
(F_{p_n}\cap U_n)$ is nowhere dense in $F_{p_n}\cap U_n$. Here, $F_{p_n}$
denotes for each $n$, the fiber over $p_n$. We can therefore
choose  $g_n\in (F_{p_n}\cap U_n)\setminus {\cal E}$ such that
$g_n(p_0)=p_n$. This yields a sequence $\{g_n\}\in U_S\setminus {\cal E}$
converging to $g$ as required.
\end{proof}

\section{Hypersurfaces complementary to $C(\gamma)$}

Our aim here is to show the existence of a complex $B$-invariant hypersurface
in the
complement of $C(\gamma)$ by constructing a function $f\in \Gamma(S,{\cal
  O}(*Y))$ and showing that the polar set of the trace transform ${\cal
  P}(Tr(f))$, lies outside of the cycle space. We continue on in the
setup and with the notation of the previous chapter.
\begin{prop}\label{trace}
If $g\in I_Y$ and there is a sequence $\{g_n\}\subset G^{\mathbb C}\setminus
I_Y$ with $g_n\rightarrow g$ and $p_n\subset g_n(C_0)\cap {\cal O}$ such
that $p_n\rightarrow p\in Y$, then there exists $f\in \Gamma(S,{\cal O}(*Y))$
with $p\in {\cal P}(f)$ and $g\in {\cal P}(Tr(f))$.
\end{prop}
\begin{proof}
Through the map
\begin{equation*}
\pi:C(\gamma)\to C(\gamma)/G^{\mathbb C}_{C_0}, \text { given by
} g\mapsto \pi(g)=C:=g(\text {c}\ell(\kappa)),
\end{equation*}
we identify the sequence $\{g_n\}\in G^{\mathbb C}\setminus I_Y$ converging to
$g\in I_Y$ with a sequence of cycles $C_n:=g_n(C_0)\not\in
I_Y/G^{\mathbb C}_{C_0}$ converging to $g(C_0)\in I_Y/G^{\mathbb C}_{C_0}$
with  $\{p_n\}\subset g_n(C_0)\cap S$ converging to $p\in Y$.

\m
Now embed $S$ (see Prop.\,\ref{embedding}) into some projective space so that
$Y$  is the intersection of $S$ with the hyperplane at infinity.
Since $g_n(C_0)\cap S\subset {\cal O}$ and this intersection is finite, we
have that it is equal to a sequence of the form $\{p_1^n,\ldots
p_{k_n}^n\}\subset {\cal O}$ with at least some $p_i^n\to p\in Y$, since
$p_n\to p\in Y$. Without loss of generality, assume $p_n=p_1^n\to\infty$. This
is equivalent to $p_1^n\to \infty $ now considered as a sequence
in $\mathbb C^m$, the complement of the hyperplane $\{z_0=0\}$ (see
Cor.\ref{hyperplane at infinity}).

\m
Since the sequences $\{p^n_i\}$ above  may be replaced by subsequences, we may
assume that $k_n=k$ is a constant.
Let $(z_1,\ldots,z_m)$ denote  standard coordinates in $\mathbb C^m$.
Then with respect to this coordinate system, $p_i^n=(z^n_{i1},\ldots
z^n_{im})$ and since $p_1^n\to \infty$, we may assume again without loss of
generality that $z^n_{11}\to \infty $ as a sequence of complex numbers.

\m
Let $a^n :=(a^n_1,\ldots,a^n_k)\in \mathbb C^k$ be the sequence of
first coordinates of $p_i^n$, in other
words, $a^n_1=z^n_{11},\ldots,z^n_{k1}$. With this translation, the
goal is to find a polynomial $P=P(z)$ of one variable so that its restriction
to ${\cal O}$ when regarded as a function $f(z):=P(z_1)$ on $\mathbb C^m$
satisfies  

\begin {equation*}
\underset{n\to \infty}{lim}\,\underset{i}{\sum}P(a^n_i)=\infty.
\end {equation*}

\m
Let ${\mathfrak S}_k$ be the symmetric group acting on $\mathbb C^k$ as
usual. Since the natural projection $\pi :\mathbb C^k\to Sym_k(\mathbb
C):=\mathbb C^k/{\mathfrak S}_k$ is proper,
it follows that $b^n:=\pi(a^n)=(a_1^n,\ldots a_k^n)$ (still denoted by
the same coordinates), diverges  in $Sym_k(\mathbb C)$. Now since the
algebraic variety $Sym_k(\mathbb C)$ is affine, it follows that there is a
regular function $R\in {\cal O}_{alg}(Sym_k(\mathbb C))$ with $R(b^n)\to
\infty $. But the regular functions on $Sym_k(\mathbb C)$ are generated
by the Newton polynomials

\begin {equation*}
N_\alpha :=\underset{i}{\sum }w_i^\alpha ,\ \alpha =1,\ldots ,k,
\end {equation*}
and consequently, $R$ is a polynomial $R=Q(N_0,\ldots ,N_k)$ in the Newton
polynomials. Since $R(b^n)\to \infty $, it follows that
at least one of the Newton polynomials, $N_p$, must be unbounded along the
sequence $\{b^n\}$. 
\m

Consequently, after going to a subsequence,

\begin {equation*}
\underset{n\to \infty}{lim}\, N_p(b^n)=\underset{n\to
  \infty}{lim}\underset{i}{\sum} (z_i^n)^p=\infty. 
\end {equation*}

Therefore the function $f(z)=P(z_1,\ldots, z_m)=z_1^p$ on $\mathbb C^m$ is such
that $\underset{n\to \infty}{lim}({\cal T}_{\psi}(f))=\infty$ and so $f$ is
the required function.
\end{proof}

The following is a converse of Prop.\,\ref{trace}.

\begin{prop}\label{polar set}
If $g\in {\cal P}({\cal T}(f))$ and the sequence $\{g_n\}\subset G^{\mathbb
  C}\setminus I_Y$ converges to $g$, then there exists a sequence
  $\{p_n\}\subset g_n(C_0)\cap S$ such that $p_n\rightarrow p\in Y$.
\end{prop}
\begin{proof}
We recall that ${\cal T}(f)$ is defined by averaging $f$ over the intersection
$g(C_0)\cap S$. Therefore, If $g_n\to g\in {\cal P}({\cal T}(f))$, then it
follows that ${\cal T}((f)(g_n(C_0)))$ also converges to ${\cal T}(f)g(C_0)$
which is infinite since  $g\in {\cal P}({\cal T}(f))$. If all the elements of
the sets $g_n(C_0)\cap S$ were bounded away from $Y$, then ${\cal T}(f)$ would
be bounded. Hence there is a point $p\in Y$ and a sequence $\{p_n\}\in
g_n(C_0)\cap S$ such that $p_n\to p$.
\end{proof}

\m
Now for $(\gamma,\kappa)$ a dual pair and $p_0\in \gamma \,\cap \,\kappa$ a
base point, define the following subset of $U_S$;
$$
U_Y:=\{g\in G^\mathbb C:g(p_0)\in Y\}.
$$

\begin{prop}\label{g in the polar set}
Let $g\in U_Y$, then there exists
$f\in\Gamma(S,{\cal O}(*Y))$ such that $g\in {\cal P}({\cal T}(f))$ and $g\not
\in C(\gamma)$.
\end{prop}
\begin{proof}
Since $g(p_0)=p\in Y$ and $D_S$ is nowhere dense in $S$, there is a sequence
of points $\{p_n\}\subset S\setminus D_S$ with $p_n$ converging to $p$. Now by
Prop.\,\ref{sequences}, there exists a sequence of transformations
$\{g_n\}\subset U_S\setminus {\cal E}$ with $g_n$ converging to $g$ such that
$g_n(p_0)=p_n$. Thus the first statement follows from Prop.$\,$\ref{trace}.

\m
Since $g\in {\cal P}({\cal T}_{\psi}(f))$, it follows that $g_n(\text
{c}\ell(\kappa))\cap S$ is a finite set say, $\{p_1^n,\ldots, p_k^n\}$ and is
contained in ${\cal O}\cap \gamma$. Since $g_n$ converges to $g\in {\cal
  P}({\cal T}_{\psi}(f))$, at least one of
the sequences $\{p_i^n\}$ converges to a point in $Y$. Assume that $p_1^n\to
p\in Y$. It now follows from Cor.$\,$\ref{sequences in S} that $g$ is not
contained in $C(\gamma)$.
\end{proof}

\m
Observe that since $I_Y$ is $K^\mathbb C$-invariant, it follows that if
$g(p_o)\in I_Y$, then there exists $k\in K^\mathbb C$ such that $g(k(p_0))\in
U_Y$. Consequently, $g\in U_Y.K^\mathbb C$.

\begin{cor}
If $g\in U_Y.K^\mathbb C$, then $I_Y$ is locally $1$-codimensional at $g$ and
is contained in the complement of $C(\gamma)$.
\end{cor}

\begin{lem}
The set $U_Y.K^\mathbb C$ is dense in $I_Y$.
\end{lem}
\begin{proof}
For any given $g\in I_Y$, there exists an arbitrarily small $h\in G^\mathbb C$
so that $hg\in U_Y.K^\mathbb C$. In particular $U_Y.K^\mathbb C\cap I_Y$ is
dense in $I_Y$.
\end{proof}

\begin{cor}
The set $I_Y$ is $1$-codimensional and contained in the complement of
$C(\gamma)$.
\end{cor}
\begin{proof}
We have shown above that if $g\in U_Y$, then $g\not \in
C(\gamma)$. Therefore, $g\not \in C(\gamma)$ if $g\in U_Y.K^\mathbb
C$. Furthermore, for any $g\in U_Y.K^\mathbb C$ there
exists $f\in\Gamma(S,{\cal O}(*Y))$ such that $g\in {\cal P}({\cal
  T}(f))$. Consequently, $U_Y.K^\mathbb C\cap I_Y$ is contained in the
complement of $C(\gamma)$ and is $1$-codimensional at each point of
$I_Y$. Since  $U_Y.K^\mathbb C\cap I_Y$ is dense in $I_Y$, it follows
that $I_Y$ is $1$-dimensional and is contained in the complement of
$C(\gamma)$.
\end{proof}

\section{Cycle spaces of nonclosed orbits}

Let $B\subset G^\mathbb C$ denote as usual an Iwasawa-Borel subgroup. Then $B$
has only finitely many orbits in $\Omega:=G^\mathbb C/K^\mathbb C$. If
$x_0:=1K^\mathbb C$ denote the base point in $\Omega$, then the orbit $B.x_0$
is open in $\Omega$, since the complexification of the Iwasawa decomposition
$G=K.A.N$ of $G$ is open in $G^\mathbb C$. The complement of the open
$B$-orbit in $\Omega$ therefore consists of a union of a finite number of
$B$-invariant complex hypersurfaces. Let $H$ be such a hypersurface invariant
for some fixed $B$, then the family $\{gH\}_{g\in G}$ consists of
$G$-translates of $H$. For some Iwasawa decomposition of $G$, the hypersurface
$H$ is  $AN$-invariant since $AN\subset B$ and consequently, this family is
equivalent to the family $\{kH\}_{k\in K}$. We denote by $\Omega_H$ the
$G$-invariant domain defined by $H$ in the following way;
$$
\Omega_H:=(\Omega\setminus \underset {k\in K}\bigcup (kH))^0,
$$
the connected component containing the base point $x_0$ in $\Omega$.

\m
We will first handle a certain situation which is present in all
nonhermitian cases and many Hermitian cases.  For this we recall the
notation in the Hermitian setting.

\m
Associated to the compact symmetric spaces $X_+$ and $X_-$ are bounded
symmetric domains realized as $G$-orbits of the neutral points in the
following way. Let $x_+$ be the neutral point in $X_+$ and $x_-$ the neutral
point in $X_-$  i.e.,
$x_+=e.P_+\in X_+$ and $x_-=eP_-\in X_-$. Let us denote by $\cal {B}$ the
bounded symmetric space $G/K$ with the
complex structure of $G.x_-$ and $\cal \bar{B}$ the bounded symmetric space
$G/K$ with the complex structure of $G.x_+$.

\m

As usual, let $(\gamma, \kappa)\in Orb_Z(G)\times Orb_Z(G^\mathbb C)$ be a
dual pair considered here in the Hermitian case for an
arbitrary flag manifold $Z$. 
Recall that by the $K^\mathbb C$-invariance of the base cycle
$C_0=\text {c}\ell(\kappa)$, the cycle
space $C(\gamma)$ is right $K^\mathbb C$-invariant.
If $C_0$ is only $K^\mathbb C$-invariant,
we regard the cycle space as
$C(\gamma )/K^\mathbb C\subset G^\mathbb C/K^\mathbb C$.
If $C_0$ is either $P_+$- or $P_-$- invariant, then we regard the cycle
space as $C(\gamma )/P_+$ or $C(\gamma)P_-$.

\m
We also recall that the Iwasawa-Borel subgroup $B$ acts on the symmetric spaces
$X_+=G^\mathbb C/P_+$ and $X_-=G^\mathbb C/P_-$. Since $b_2(X)=1$, there is a
unique $B$-invariant hypersurface $H_0$ in the complement of $\cal {B}$.
We will maintain notation of the previous sections and recall that we have
the $B$-invariant hypersurface $I_Y$ containing ${\cal P}({\cal
  T}_{\psi}(f))$, the polar 
set of the trace transform ${\cal T}_{\psi}(f)$, constructed in Section 7.
Moreover, $I_Y$ is in the complement of $C(\gamma)$.

\m
Now define $H$ to be the maximal $B$-invariant hypersurface in the complement
of the cycle space $C(\gamma)$. In the proof of the main results, it will be
important if the hypersurface $H$
is a $\pi_+$- (or $\pi_-$)-lift of $H_0$ from $X_+$ (or $X_-$) or not and so
these two cases will be distinguished.

\m

Since the $G^\mathbb C$-isotropy subgroup $G^{\mathbb C}_{C_0}=\{g\in G^\mathbb
C:g(C_0)=C_0\}$ of the base cycle $C_0$ contains $K^\mathbb C$, the orbit
$G^\mathbb C.C_0$ is either $G^\mathbb C/\tilde {K}^\mathbb C$, where $\tilde
{K}^\mathbb C$ is a finite extension of $K^\mathbb C$ or it is one of the
compact Hermitian symmetric spaces $G^\mathbb C/P_+$ or $G^\mathbb C/P_+$.

\subsection {Case {\rm I}: $H$ is not a lift}

Here we consider the case when there exists a maximal $B$-invariant
hypersurface $H$ in the
complement of $C(\gamma)$ which is not a lift. This means that
for every Iwasawa-Borel group $B$ the maximal $B$-invariant 
hypersurface is neither of the form $H=\pi ^{-1}_+(H_+)$ nor of the
from $H=\pi ^{-1}_-(H_-)$, where $\pi _\pm: G^\mathbb C/K^\mathbb C\to
G^\mathbb C/P_\pm =X_\pm $ are the standard projections.  Here
$H_+$ (resp. $H_-$) is the unique $B$-invariant hypersurface in
$X_+$ (resp. $X_-$).

\s
Of course in the nonhermitian case there are no such projections
and therefore this imposes no condition.

\m
Since $H$ is not a lift and $\Omega _H$ contains
the cycle space, the orbit $G^\mathbb C.C_0$ in ${\cal C}_q(Z)$ is $G^\mathbb
C/\tilde {K}^\mathbb C$. Making use of the finite covering map $\pi:G^\mathbb
C/K^\mathbb C \to G^\mathbb C/\tilde {K}^\mathbb C$, we lift the cycle space
to $C(\gamma)/K^\mathbb C$ which is an open subset of $G^\mathbb C/K^\mathbb
C$ and still denote it by $C(\gamma)$. 

\m
This will allow us to be able to compare $C(\gamma)$ with the domain
$\Omega_{AG}$ which is contained in $G^\mathbb C/K^\mathbb C$. We state the
following result in this context.

\begin {thm} \label {hypersurface is not a lift}
If the maximal $B$-invariant hypersurface $H$ in the complement of $C(\gamma)$
is not a lift, then
$$
\Omega _{AG}=C(\gamma)=\Omega_H.
$$
\end {thm}

In order to prove this theorem, we first state some known results concerning
an open $G$-orbit $\gamma_{open}\in Orb_Z(G)$.

Firstly, it has been shown in \cite {GM} that the intersection of all cycle
spaces $C(\gamma)$ as $\gamma$ ranges over $ Orb_Z(G)$ and $Q$ ranges over all
parabolic subgroups of $G^\mathbb C$ is the same as the intersection of all
cycle spaces $C(\gamma_{open})$ for all open $G$-orbits in $Z=G^\mathbb C/B$,
for $B$ a Borel subgroup of $G^\mathbb C$. That is

\begin{lem} (\cite {GM}).\label {Lemma 9.2}
$$
\underset {\underset {Z=G^\mathbb C/Q}{\gamma \in Orb_Z(G)}}\bigcap
C(\gamma)=\underset {\underset {Z=G^\mathbb C/B}{\gamma \text {open}}}\bigcap
C(\gamma).
$$
\end{lem}

If $G$ is of Hermitian type the following result is also known

\begin{thm} (\cite {HW1},\cite {WZ}).
If $G$ is of Hermitian type and $\gamma \in Orb_Z(G)$ is open, then either
\begin{enumerate}
\item {the base cycle $\kappa$ is $P_+$-invariant
(resp. $P_-$-invariant)\\
and $C(\gamma)=\bar {\cal B}$ (resp. ${\cal B}$)}\\
or
\item {the base cycle is only invariant by $K^\mathbb C$
  and $C(\gamma)={\cal B}\times{\cal {\bar B}}$}.
\end{enumerate}
\end{thm}

Since $\Omega _{AG}={\cal B}\times \bar {\cal B}$ in the Hermitian
case (\cite {BHH}), we have the following consequence.

\begin{cor}
If $G$ is of Hermitian type and $\gamma \in Orb_Z(G)$ is open, then
$$
\underset {\underset {Z=G^\mathbb C/B}{\gamma \text {open}}}\bigcap C(\gamma)
={\cal B}\times{\cal {\bar B}}= \Omega_{AG}.
$$
\end{cor}
\begin {proof}
Although there will always be open orbits with cycle spaces of
the second type above, it suffices to note that

\begin {equation*}
\pi ^{-1}_-({\cal B})\,\cap \pi ^{-1}_+\bar {\cal B}=
{\cal B}\times  \bar {\cal B}
\end {equation*}

which is embedded as $\Omega _{AG}$ in $G^\mathbb C/K^\mathbb C$.
\end {proof}

If $G$ is not of Hermitian type, then we have the following result
(\cite {FH}) for open $G$-orbits in $Z=G^\mathbb C/B$.

\begin{thm} (\cite {FH}).
If $G$ is not Hermitian, then 
$$C(\gamma)=\Omega_{AG}
$$
for every open $G$-orbit $\gamma \in Orb_Z(G)$ and $Z$ any flag manifold.
\end{thm}

\begin{cor}
In both the Hermitian and non-Hermitian cases,
$$
\underset {\underset {Z=G^\mathbb C/B}{\gamma \text {open}}}\bigcap C(\gamma)
=\Omega_{AG}.
$$
\end{cor}

\begin{cor}\label {U contained in cycle space}
For all $\gamma \in Orb_Z(G)$,
$$
\Omega_{AG}\subset C(\gamma).
$$
\end{cor}
\begin{proof}
Apply Lemma \ref {Lemma 9.2}.
\end{proof}

This immediately implies that if $H$ is a maximal $B$-invariant hypersurface
which lies outside $C(\gamma)$, then
$$
\Omega_{AG}\subset C(\gamma)\subset \Omega_H
$$.

\m
{\it Proof of Theorem \ref{hypersurface is not a lift}}.   
This is now an immediate consequence of the work in \cite {FH}. There as a
first step it is shown that if $H$ is not a lift from either $G^\mathbb C/P_+$
or $G^\mathbb C/P_-$, then the domain $\Omega_H$ is Kobayashi hyperbolic. The
following main result of \cite {FH} then completes our proof: If $\hat \Omega$
is a $G$-invariant, Kobayashi hyperbolic Stein domain in $\Omega$ which
contains $\Omega_{AG}$, then $\hat \Omega =\Omega_{AG}$.
\hfill $\square$

\b
Theorem \ref{hypersurface is not a lift} proves
our main result Theorm \ref{main theorem} for non-closed orbits in the case
where tha base cycle $C_0=\text {c}\ell(\kappa)$ is neither $P_-$- nor
$P_-$-invariant .

\m
The following is a consequence of Theorem \ref{hypersurface is not a lift}.

\begin {cor} \label{non Hermitian case}
If $G$ is not of Hermitian type, then
\begin {equation*}
C(\gamma )=\Omega _{AG}
\end {equation*}
for all $\gamma \in Orb_Z(G)$.
\end {cor}

\begin {proof}
Since $G$ is non-Hermitian, no $H$ is a lift.
\end {proof}

\subsection {Case {\rm II}: Every $H$ is a lift}

We now consider the case where the maximal $B$-invariant hypersurface
$H=\pi^{-1}(H_+)$ is a lift from $G^\mathbb C/P_+$, where $\pi _+$ is the
natural projection. Of course the discussion is the same if every $H$ is a
lift from $G^\mathbb C/P_-$. 

\m
We begin by proving the following

\begin {thm}
If $\text {c}\ell(\kappa)$ is not $P_+$-invariant, then no Schubert variety
$S={\cal O}\, \dot \cup \, Y \in {\cal S}_\kappa$ defines a maximal
$B$-invariant hypersurface $H$ which is a lift from $G^\mathbb C/P_+$.
\end {thm}
\begin {proof}
Suppose to the contrary that there is some Schubert variety $S={\cal O}\,\dot
\cup \, Y\in {\cal S}_\kappa$ defining $H=\pi^{-1}(H_+)$ which is lift. This
is equivalent to the domain $\Omega_H =\pi ^{-1}_+(\Omega_{H_+})$ being a
lift.

\m
Let $x_0\in \gamma$ be the base point with $\kappa=K^\mathbb
C.x_0$. Since $\kappa$ is not $P_+$-invariant, $\text {c}\ell(P_+.x_0)$
contains $\text {c}\ell(\kappa)$ as a proper subvariety. Now the intersection
$\text {c}\ell(\kappa)\cap S\subset {\cal O}$ and is transversal in $Z$. Thus
every component of $P_+.x_0\cap {\cal O}$ is positive dimensional. Since
${\cal O}=\mathbb C^{m({\cal O})}$ is affine, every such component has at
least one point of $Y$ in its closure.

\m
Thus for every arbitrarily small neigborhood $U$ of the identity in $G^{\mathbb
  C}$ there exists $h\in P_+$, and $g\in U$ with $gh.x_0\in Y$.  Consequently,
 $gh\in U_S$ and it follows from Prop.$\,$\ref {g in the polar set} that
$ghK^\mathbb C$ is in the maximal $B$-invariant hypersurface $H$. Hence
$ghK^\mathbb C\not \in \Omega_H =\pi ^{-1}_+(\Omega_{H_+})$.  

\m
On the other hand $C(\gamma )\subset \Omega _H$ and $C(\gamma )$
contains an open neighborhood $U$ of the identity.  Consequently,
$ghK^\mathbb C\subset \Omega _H$ for every $g\in U$ and $h\in P_+$.
Thus we have reached a contradiction, and therefore no Schubert variety $S\in
{\cal S}_\kappa$ defines $H=\pi_+^{-1}(H_+)$.
\end {proof}

\m
It therefore follows that if $\text {c}\ell(\kappa)$ is neither $P_+$- nor
$P_-$-invariant,
the domain $\Omega_H$ is Kobayashi hyperbolic \cite {FH}. Thus the proof of
our main result Theorem \ref{main theorem} for
non-closed orbits in this case is completed just like in the proof of Theorem
\ref{hypersurface is not a lift}.

\m
By taking contrapositions in the above theorem, we obtain

\begin {cor} \label {hypersurface is a lift}
If the maximal $B$-invariant hypersurface $H=\pi_+^{-1}(H_+)$ is a lift, then
$\text {c}\ell(\kappa)$ is $P_+$-invariant.
\end {cor}

\begin{cor}
If the maximal $B$-invariant hypersurface $H=\pi_+^{-1}(H_+)$ (resp. $H=\pi
^{-1}_-(H_-)$) is a lift from $G^\mathbb C/P_+$ (resp. $G^\mathbb C/P_-$),
then $C(\gamma)={\cal {\bar B}}$ (resp. $C(\gamma )={\cal B}$).
\end {cor}
\begin{proof}
We have seen above that under this assumption 
$\text {c}\ell(\kappa)$ is $P_+$-invariant. This
implies that $g(\text {c}\ell(\kappa))$ is
$gP_+g^{-1}$-invariant. Consequently, if $g\in C(\gamma)$, then
$gP_+g^{-1}.g\subset C(\gamma)$, that is, $C(\gamma)$ is right
$P$-invariant. Thus $C(\gamma)$ may be regarded as a domain in $G^\mathbb
C/P_+$. Since ${\cal {\bar B}}\subset G^\mathbb C/P_+$ is a 
$G$-orbit, it follows that  ${\cal {\bar B}}\subset C(\gamma)$. Since
$C(\gamma)\subset \Omega_{H_+}$ and we know that $\Omega_{H_+}={\cal {\bar B}}$
(see for example \cite {H}), the result follows.
\end{proof}

\m
This completes the proof of our main result Theorem \ref{main theorem} for
non-closed orbit.

\m
As a consequence of the work in this and the previous subsections we
now have the following result.

\begin {thm}
Suppose $\gamma $ is a nonclosed $G$-orbit. If $G$ is of Hermitian type and
$\text {c}\ell(\kappa)$ is neither $P_+$- nor $P_-$-invariant, then
\begin {equation} 
C(\gamma )=\Omega _{AG}.
\end {equation}
If $\gamma $ is nonclosed and $\text{c}\ell (\kappa )$ is
$P_+$-invariant (resp. $P_-$-invariant), then
$C(\gamma )=\bar {\cal B}$ (resp. $C(\gamma )={\cal B}$).
\end {thm}

Note that if $\text{c}\ell (\kappa )$ is $P_+$- or $P_-$-invariant,
then the orbit $G^\mathbb C.C_0$ in the cycle space ${\cal C}_q(Z)$ is
$G^\mathbb C/P_+$ or $G^\mathbb C/P_-$.  Thus the latter statement,
$C(\gamma )=\bar {\cal B}$ or $C(\gamma )={\cal B}$, is a statement
in ${\cal C}_q(Z)$.

\m
The former statement must be interpreted.  In that case
$G^\mathbb C.C_0=G^\mathbb C/\tilde K^\mathbb C$, where
$\tilde K^\mathbb C$ is possibly a finite extension of $K^\mathbb C$.
If we regard the sets $\Omega _H$, which are defined by incidence
geometry, as being in ${\cal C}_q(Z)$, then the cycle space statement
is $C(\gamma )=\Omega _H$ in ${\cal C}_q(Z)$.

\s
However, by the main result of (\cite {FH}) the lift of $\Omega _H$
in $G^\mathbb C/K^\mathbb C$ is $\Omega _{AG}$, and, since 
$\Omega _{AG}$ is a cell, this lift is biholomorphic.  Thus
in this sense we write $C(\gamma )=\Omega _{AG}$ in $G^\mathbb C/K^\mathbb C$.

\section {Cycle spaces of closed orbits}

Here we consider the case of the closed $G$-orbit $\gamma _{\text{c}\ell }$
and its dual $\kappa _{op}$ which is the open $K^\mathbb C$-orbit 
in $Z$.  Duality is just the statement that
$\kappa _{op}\supset \gamma _{\text{c}\ell }$.

\b
We begin by recalling the behavior of duality with respect to the
partial ordering of orbits defined by the closure operation.  For
this, if ${\cal O}_1$ and ${\cal O}_2$ are orbits with
${\cal O}_2\subset \text{c}\ell({\cal O}_1)\setminus {\cal O}_1$,
we write ${\cal O}_1<{\cal O}_2$.  The following is a well-known
aspect of the duality theory.

\begin {prop} \label {orbit ordering}
If $(\gamma _1,\kappa _1)$ and $(\gamma _2,\kappa _2)$ are dual
pairs, then
\begin {equation*}
\gamma _1<\gamma _2\ \Leftrightarrow \ \kappa _2<\kappa _1.
\end {equation*}
\end {prop} 
{\it Sketch of proof.} Suppose that $\gamma _1<\gamma _2$.
As a consequence, $\gamma _1\,\cap \, \kappa _2\not=\emptyset $.
Recall that $\kappa _1\,\cap \, \gamma _1$ is realized as a
strong deformation retract by the gradient flow
$\varphi _t$ of the norm of a certain moment map
(see, e.g., [BL]).  This flow is $K$-invariant and tangent to
all $K^\mathbb C$- and $G$-orbits.

\m
Now take $p\in \kappa _2\, \cap \, \gamma _1$ and let
$q:=\underset{t\to \infty}{\lim }\varphi _t(p)
\in \kappa _1\cap \gamma _1$. Since $\text{c}\ell (\kappa _2)$
is invariant under the flow, it follows that
$\kappa _1=K^\mathbb C.q\subset \text{c}\ell(\kappa _2)$.  By assumption
$\gamma _1\not=\gamma _2$.  Thus $\kappa _1\subset \text{c}\ell (\kappa
_2)\setminus \kappa _2$ as required.

\s
The converse implication goes in exactly the same way.
\hfill $\square $

\m
This result will now be used in the following special situation.
Let
\begin {equation*}
bd(\kappa _{op})=Z\setminus \kappa _{op}=
A_1\, \cup \ldots \cup \, A_k
\end {equation*}
be the decomposition of the boundary $bd(\kappa _{op})$ of $\kappa_{op}$ as a
union of its irreducible components. In each $A_j$ there
is a unique, Zariski dense open $K^\mathbb C$-orbit $\kappa _j$ with
dual $G$-orbit $\gamma _j$, $j=1,\ldots ,k$.

\begin {cor}
For every $j$ it follows that 
$\text{c}\ell (\gamma _j)=\gamma _j\,\dot \cup \, \gamma _{\text{c}\ell }$.
\end {cor}
\begin {proof}
If $\gamma $ is a $G$-orbit contained in
$\text{c}\ell (\gamma _j)\setminus \gamma _j$, then its dual
$K^\mathbb C$-orbit $\kappa $ satisfies $\kappa >\kappa _j$.  But
$\kappa =\kappa _{\text{op}}$ is the only $K^\mathbb C$-orbit with this
property.
\end {proof}

\begin {cor}
For $\gamma _j$ as above,
\begin {equation*}
C(\gamma _{\text{c}\ell })\subset C(\gamma _j).
\end {equation*}
\end {cor}
\begin {proof}
If $g\in \text{bd}(C(\gamma _j))$, then
$g(\text{c}\ell (\kappa _j))\, \cap \, \text{bd}(\gamma _j)\not =\emptyset $.
Thus $g(\text{c}\ell (\kappa _j)) \, \cap \, 
\gamma _{\text {c}\ell }\not =\emptyset $.  In particular,
$g\not \in C(\gamma _{\text{c}\ell })$.
\end {proof}
\begin {thm}
For every $Z=G/Q$ it follows that
\begin {equation*}
C(\gamma _{\text{c}\ell })=\Omega _{AG}.
\end {equation*}
\end {thm}
\begin {proof}
By the above Corollary and Cor.$\,$\ref{U contained in cycle space}, for
all $j$
\begin {equation*}
\Omega _{AG}\subset C(\gamma _{\text{c}\ell })\subset C(\gamma _j).
\end {equation*}

\m
If $C(\gamma _j)=\Omega _{AG}$ for some $j$, then the proof is
finished.  

\s
In the Hermitian case, if, e.g.,
$C(\gamma _1)=\pi ^{-1}_-({\cal B})$ and 
$C(\gamma _2)=\pi ^{-1}_+(\bar {\cal B})$, then
$C(\gamma _1)\,\cap \,C(\gamma _2)=\Omega _{AG}$ and the
proof is finished in that case as well.

\m
Thus we may assume that $C(\gamma _j)=\pi _-^{-1}({\cal B})$
for all $j$, or equivalently that $\text{c}\ell (\kappa _j)$
is $P_-$-invariant for all $j$.  But this is in turn equivalent
to $\kappa _{\text{op}}$ being $P_-$-invariant.

\m
However, $\kappa _{\text{op}}$ can not be $P_-$-invariant.  To see
this, note that if it were invariant, then $C(\gamma _{\text{c}\ell })$
would be $P_-$-invariant.  Since
$\Omega _{AG}\subset C(\gamma _{\text{c}\ell })$, 
this would imply that
$P_-.\Omega _{AG}\subset C(\gamma _{\text{c}\ell })$.

\s
But the $P_-$-orbit of a generic point in $\Omega _{AG}$ is Zariski
open in $\Omega $.  Consequently, if $\kappa _{\text{op}}$ were
$P_-$-invariant, it would follow that $C(\gamma _{\text {c}\ell })$
would contain a Zariski open subset of $\Omega $.  By the identity
principle, this is contrary to the complement of $C(\gamma _{\text{c}\ell })$
being nonempty and $G$-invariant.
\end {proof}

\m
This completes the proof of our main theorem, Theorem \ref{main theorem}.

\end{document}